\documentclass[a4,11pt]{amsart}

\numberwithin{equation}{section}

\newtheorem{mainthm}{Theorem}
\newtheorem{mainprop}[mainthm]{Proposition}

\newtheorem{lem}{Lemma}[section]
\newtheorem{thm}[lem]{Theorem}
\newtheorem{prop}[lem]{Proposition}

\newtheorem{rem}[lem]{\it Remark}

\newcommand{\conj}{\overline}
\newcommand{\cz}{{\conj{z}}}
\newcommand{\cw}{{\conj{w}}}
\newcommand{\cv}{{\conj{v}}}
\newcommand{\pa}{\partial}
\newcommand{\wt}{\widetilde}

\renewcommand{\th}{\theta}
\newcommand{\e}{\varepsilon}
\newcommand{\im}{i}
\newcommand{\qtext}[1]{\quad\text{#1}\ \ }

\newcommand{\zs}{\setminus\{0\}}

\newcommand{\C}{\mathbb{C}}
\newcommand{\R}{\mathbb{R}}

\newcommand{\N}{\mathbb{N}}
\newcommand{\calA}{\mathcal{A}}
\newcommand{\calB}{\mathcal{B}}
\newcommand{\calC}{\mathcal{C}}
\newcommand{\calE}{\mathcal{E}}

\newcommand{\calO}{\mathcal{O}}

\renewcommand{\Im}{\operatorname{Im}}
\newcommand{\Vol}{\operatorname{Vol}}

\begin{document}

\title[ Logarithmic singularity of the Szeg\"o kernel
]{
 Logarithmic singularity of the Szeg\"o kernel and 
a global invariant of strictly pseudoconvex domains
 }

\author[]{Kengo Hirachi}

\address{
Graduate School of Mathematical Sciences\\
The University of Tokyo\\
Komaba, Megro, Tokyo 153-8914, Japan}
\email{hirachi@ms.u-tokyo.ac.jp}

%


\maketitle

\section{Introduction} 

This paper is a continuation of Fefferman's program \cite{F2} for studying the
geometry and analysis of strictly pseudoconvex domains. The key idea of the program is to consider the Bergman and Szeg\"o kernels of the domains as analogs of the heat kernel of Riemannian manifolds. In Riemannian (or conformal) geometry, the coefficients of the asymptotic expansion of the heat kernel can be expressed in terms of the curvature of the metric; by integrating the coefficients one obtains index theorems in various settings. For the Bergman and Szeg\"o kernels, there has been much progress made on the description of their asymptotic expansions based on invariant theory (\cite{F2}, \cite{BEG}, \cite{Hi3}); we now seek for invariants that arise
from the integral of the coefficients of the expansions.

We here prove that the integral of the coefficient of the logarithmic singularity of the Szeg\"o kernel gives a biholomorphic invariant of a domain $\Omega$, or a CR invariant of the boundary $\pa\Omega$, and moreover that the invariant is unchanged under perturbations of the domain (Theorem \ref{thmSzego}). We also show that the same invariant appears as the coefficient of the logarithmic term of the volume expansion of the domain with respect to the Bergman volume element (Theorem \ref{thmBergman}). This second result is an analogy of the derivation of a conformal invariant from the volume expansion of conformally compact Einstein metrics which arises in the AdS/CFT correspondence -- see \cite{Gr4} for a discussion and references.

The proofs of these results based on Kashiwara's microlocal analysis of the Bergman kernel in \cite{Kas}, where he showed that the reproducing property of the Bergman kernel on holomorphic functions can be ``quantized" to a reproducing property of the microdifferential operators (i.e., classical analytic pseudodifferential operators).
It provides a system of microdifferential equations that characterizes the singularity of the Bergman kernel (which can be formulated as a microfunction) up to a constant multiple; such an argument can be equally applied to the Szeg\"o kernel. These systems of equations are used to overcome one of the main difficulties, when we consider the analogy to the heat kernel, that the Bergman and Szeg\"o kernels are not defined as solutions to differential equations.

Let $\Omega$ be a relatively compact, smoothly bounded strictly pseudoconvex domain in a complex manifold $M$. We take a pseudohermitian structure $\th$, or a contact form, of $\pa\Omega$ and define a surface element $d\sigma=\th\wedge (d\th)^{n-1}$. Then we may define the Hardy space $\calA(\pa\Omega,d\sigma)$ consisting of the  boundary values of holomorphic functions on $\Omega$ that are $L^2$ in the norm $\|f\|^2=\int_{\pa\Omega}|f|^2 d\sigma$.  The Szeg\"o kernel $S_\th(z,\conj{w})$ is defined as the reproducing kernel of $\calA(\pa\Omega,d\sigma)$, which can be extended to a holomorphic function of $(z,\conj{w})\in\Omega\times\conj\Omega$ and has a singularity along the boundary diagonal. If we take a smooth defining function $\rho$ of the domain, which is positive in $\Omega$ and $d\rho\ne0$ on $\pa\Omega$,
then (by \cite{F0} and \cite{BS}) we can expand the singularity as
\begin{equation}\label{S-exp}
 S_\th(z,\conj z)=\varphi_\th(z)\rho(z)^{-n}+\psi_\th(z)\log\rho(z),
\end{equation}
where $\varphi_\th$ and $\psi_\th$ are functions on $\Omega$ that are smooth up to the boundary. Note that $\psi_\th|_{\pa\Omega}$ is independent of the choice of $\rho$ and is shown to gives a local invariant of the pseudohermitian structure $\th$.

\begin{mainthm} \label{thmSzego}
{\rm (i)}
The integral 
$$
  L(\pa\Omega,\th)=\int_{\pa\Omega}\psi_\th\,\th\wedge (d\th)^{n-1}
$$
is independent of the choice of a pseudohermitian structure $\th$ of $\pa\Omega$. Thus we may write $L(\pa\Omega)=L(\pa\Omega,\th)$. 

\smallskip

{\rm (ii)}
Let $\{\Omega_t\}_{t\in\R}$ be a $C^\infty$ family of strictly pseudoconvex domains in $M$.  Then $L(\pa\Omega_t)$ is independent of $t$.
\end{mainthm}

In case $n=2$, we have shown in \cite{Hi2} that 
$$
\psi_\th|_{\pa\Omega}=\frac1{24\pi^2} 
(\Delta_bR-2\Im A_{11,}{}^{11}), 
$$
where $\Delta_b$ is the sub-Laplacian, $R$ and $A_{11,}{}^{11}$ are respectively the scalar curvature and the  second covariant derivative of the torsion of the Tanaka-Webster connection for $\th$.
Thus the integrand $\psi_\th\,\th\wedge d\th$ is nontrivial and does depend on $\th$, but it also turns out that $L(\pa\Omega)=0$ by Stokes' theorem.
For higher dimensions, we can still give examples of $(\pa\Omega, \th)$ for which $\psi_{\th}|_{\pa\Omega}\not\equiv0$. However, the evaluation of the integral
is not easy and, so far, we can only give examples with trivial $L(\pa\Omega)$ -- see Proposition \ref{prop-tube-volume} below.

We were led to consider the integral of $\psi_\th$ by the works of Branson-{\O}rstead \cite{BO} and Parker-Rosenberg \cite{PR} on the constructions of conformal invariant from the heat kernel $k_t(x,y)$ of conformal Laplacian, and their CR analogue for CR invariant sub-Laplacian by Stanton \cite{Stanton}. For a conformal manifold of 
even dimension $2n$ (resp.~CR manifold of dimension $2n-1$), the integral of
the coefficient $a_n$ of the asymptotic expansion $k_t(x,x)\sim t^{-n}\sum_{j=0}^\infty a_j(x) t^j$ is shown to be a conformal (resp.\ CR) invariant,
while the integrand $a_ndv_g$ does depend on the choice of a scale $g\in[g]$
 (resp.\ a contact form $\th$).
This is a natural consequence of the variational formula for the kernel $k_t(x,y)$ under conformal scaling, which follows from the heat equation. Our Theorem \ref{thmSzego} is also a consequence of a variational formula of the Szeg\"o kernel, which is obtained as a part of system of microdifferential equations for the family of Szeg\"o kernels (Proposition \ref{Dt-prop}).

We next express $L(\pa\Omega)$ in terms of the Bergman kernel. Take a $C^\infty$ volume element $dv$ on $M$. Then the Bergman kernel $B(z,\cw)$ is defined as the reproducing kernel of the Hilbert space 
$\calA(\Omega,dv)$ of $L^2$ holomorphic functions on $\Omega$ with respect to $dv$. The volume of $\Omega$ with respect to the volume element $B(z,\conj z)dv$ is infinite. We thus set $\Omega_\e=\{z\in\Omega:\rho(z)>\e\}$ and consider the asymptotic behavior of 
$$
\Vol(\Omega_\e)=\int_{\Omega_\e}B(z,\conj z)\,dv
$$
as $\e\to+0$.

\begin{mainthm}\label{thmBergman}
 For any volume element $dv$ on $M$ 
 and any defining function $\rho$ of $\Omega$, 
 the volume $\Vol(\Omega_\e)$ admits an expansion 
\begin{equation}\label{volume-exp}
\Vol(\Omega_\e)=\sum_{j=0}^{n-1} C_j\, \e^{j-n}+
L(\pa\Omega) \log\e+O(1),
\end{equation}
where $C_j$ are constants, $L(\pa\Omega)$ is the invariant given in Theorem $1$ and $O(1)$ is a bounded term.
\end{mainthm}

The volume expansion \eqref{volume-exp}
can be compared with that of  
conformally compact Einstein manifolds (\cite{HS}, 
 \cite{Gr4}); there one considers
a complete Einstein metric  $g_+$ on the interior $\Omega$ of a
compact manifold with boundary and a conformal structure
$[g]$ on $\pa\Omega$, which is obtained as a scaling limit
of $g_+$. For each choice of a preferred defining function $\rho$
corresponding to a conformal scale,
we can consider the volume expansion of the form \eqref{volume-exp}
with respect to $g_+$. If  $\dim_\R\pa\Omega$ is even, the coefficient of the logarithmic term is shown to be a conformal invariant of the boundary $\pa\Omega$.  
Moreover, it is shown in \cite{GZ} and \cite{FG2} that this conformal invariant
can be expressed as the integral of Branson's $Q$-curvature
\cite{Br}, a local Riemannian invariant which naturally arises from conformally invariant differential operators. 
We can relate this result with ours via Fefferman's
Lorentz conformal structure defined on a $S^1$-bundle over
the CR manifold $\pa\Omega$. In case $n=2$,
we have shown in \cite{FH} that $\psi_\th|_{\pa\Omega}$ agrees with the $Q$-curvature of the Fefferman metric; while such a relation is not known for
higher dimensions.

So far, we have only considered the coefficient $L(\pa\Omega)$ of the expansion \eqref{volume-exp}. But other coefficients may have some geometric meaning if one choose $\rho$ properly; here we mention one example. Let $E\to X$ be a positive hermitian line bundle over a compact complex manifold $X$ of dimension $n-1$; then the unit tube in the dual bundle $\Omega=\{v\in E^*:  |v|<1\}$ is strictly pseudoconvex.
We take $\rho=-\log |v |^2$ as a defining function of $\Omega$ and define a volume element on $E^*$ by $dv=i^n\pa\rho\wedge\conj\pa\rho\wedge(\pa\conj\pa\rho)^{n-1}/(n-1)!$.

\begin{mainprop} \label{prop-tube-volume}
Let $B(v,\cv)$ be the Bergman kernel of $\calA(\Omega,dv)$. Then the volume of the domain $\Omega_\e=\{v\in E^*|\rho(v)>\e\}$ with respect to  the volume element $Bdv$ satisfies
\begin{equation}\label{tube-volume}
 \Vol(\Omega_\e)=2\pi
 \int_0^\infty e^{-\e\, t}P(t)dt+O(\e^\infty).
\end{equation}
Here $P(t)$ is the Hilbert polynomial of $E$, which is determined by the condition 
$P(m)=\dim H^0(M,E^{\otimes m})$ for $m\gg 0$. 
\end{mainprop}

This formula suggests a link between the expansion of $\Vol(\Omega_\e)$ and index theorems. But in this case the right-hand side of \eqref{tube-volume} does not contain $\log\e$ term and hence $L(\pa\Omega)=0$. (Note that $dv$ is singular along the zero section,
but we can modify it to a $C^\infty$ volume element without changing
\eqref{tube-volume}; see Remark \ref{rem-volume}.)

Finally, we should say again that we know no example of a domain with nontrivial $L(\pa\Omega)$ and need to ask the following:

\medskip
\noindent
{\bf Question.}{\em \/
Does there exist a strictly pseudoconvex domain $\Omega$ such that $L(\pa\Omega)\ne0?$
}
\medskip

This paper is organized as follows. In \S\ref{sec-microfunction}, we formulate the Bergman and Szeg\"o kernels as microfunctions.  We here include a quick review
of the theory of microfunctions in order that the readers can grasp the arguments of this paper even if one is unfamiliar with the subject. In \S\ref{sec-Kashiwara} we recall Kashiwara's theorem on the microlocal characterization of the Bergman and Szeg\"o kernels and derive a microdifferential relation between the two kernels and  a first variational formula of the Szeg\"o kernel. Under these preparations, we give in \S\ref{sec-proofs} the proofs of the main theorems. 
Finally in \S\ref{sec-tube}, we prove Proposition \ref{prop-tube-volume} by following Catlin \cite{Cat} and Zelditch \cite{Zel}, where the singularities of the Bergman and Szeg\"o kernels of the disk bundle are related to the asymptotic behavior of the Bergman kernels of the sections of $E^{\otimes m}$ as $m\to\infty$. With our choice of defining function $\rho$, we can simplify their arguments and give an explicit formula for this correspondence  in the course of the proof -- see \eqref{BsimA}.

\section{The Bergman and Szeg\"o kernels as microfunctions}
\label{sec-microfunction}

In this preliminary section, we explain how to formulate the theorems in terms of microfunctions, which are the main tools of this paper. We here recall all the definitions and results we use from the theory of microfunctions,
with an intention to make this section introductory to the theory. 
A fundamental reference for this section is  Sato-Kawai-Kashiwara \cite{SKK}, but a concise review of the theory by Kashiwara-Kawai \cite{KK} will be sufficient for understating the arguments of this paper. For comprehensive introductions to microfunctions and microdifferential operators, we refer to \cite{KKK}, \cite{Schapira} and \cite{Kan}.
 
\subsection{Singularity of the Bergman kernel}
We start by recalling the form of singularity of the Bergman kernel, which naturally lead us to the definition of homomorphic microfunctions.

 Let $\Omega$ be a strictly pseudoconvex domain in a complex manifold $M$ with $C^\omega$ boundary $\pa\Omega$.  We denote by $M_\R$ the underlying $C^\omega$ manifold and its complexification by $X=M\times\conj M$ with imbedding $\iota:M\to X$, $\iota(z)=(z,\cz)$.
We fix a $C^\omega$ volume element $dv$ on $M$ and define the Bergman kernel as the reproducing kernel of $\calA(\Omega,dv)=L^2(\Omega,dv)\cap \calO(\Omega)$,
where $\calO$ denotes the sheaf of holomorphic functions. Clearly we have $B(z,\cw)\in\calO(\Omega\times\conj\Omega)$, while it has singularity on the boundary diagonal.
If we take a defining function $\rho(z,\cz)$ of $\pa\Omega$, then at each boundary point $p\in\pa\Omega$, we can write the singularity of $B(z,\cw)$ as
$$
 B(z,\cw)=\varphi(z,\cw)\rho(z,\cw)^{-n-1}+\psi(z,\cw)\log\rho(z,\cw).
$$
Here $\rho(z,\cw)$ is the complexification of $\rho(z,\cz)$ and  $\varphi,\psi\in\calO_{X,p}$, where $p$ is identified with $\iota(p)\in X$. 
Moreover it is shown that this singularity is locally determined:
if $\Omega$ and $\wt\Omega$ are strictly pseudoconvex domains that agree near a boundary point $p$, then $B_\Omega(z,\cw)-B_{\wt\Omega}(z,\cw)\in \calO_{X,p}$. See \cite{Kas} and
Remark \ref{Kas-rem} below. Such an $\calO_X$ modulo class plays an essential role in the study of the system of differential equations and is called a holomorphic microfunction, which we define below in more general setting.

\subsection{Microfunctions: a quick review}
Microfunctions are ``singular part" of holomorphic functions on wedges at the edges. To formulate them, we first introduce the notion of hyperfunctions, which are generalized functions obtained by the sum of ``ideal boundary values" of holomorphic functions. 

For an open set $V\subset\R^n$ and an open convex cone
$\Gamma\subset\R^n$, we denote by $V+\im \Gamma0\subset\C^n$
an open set that asymptotically agrees with the wedge $V+\im \Gamma$
at the edge $V$. The space of {\em hyperfunctions} on $V$ is defined as a vector space of formal sums of the form
\begin{equation}\label{hyper-f}
 f(x)=\sum_{j=1}^m F_j(x+\im\Gamma_j0),
\end{equation}
where $F_j$ is a holomorphic function on $V+\im \Gamma_j0$, that allow the reduction 
$
F_j(x+\im\Gamma_j0)+F_k(x+\im\Gamma_k0)
=F_{jk}(x+\im\Gamma_{jk}0),
$
where $\Gamma_{jk}=\Gamma_j\cap\Gamma_k\neq\emptyset$ and
$F_{jk}=F_j|_{\Gamma_{jk}}+F_k|_{\Gamma_{jk}}$, and its reverse conversion. We denote the sheaf of hyperfunctions by $\calB$. Note that if each $F_j$ is of polynomial growth in $y$ at $y=0$ (i.e. $|F_j(x+iy)|\le \text{const.}|y|^{-m} $), then
$\sum_j\lim_{\Gamma_j\ni y\to 0}F_j(x+iy)$ converges to a distribution $\wt f(x)$ on $V$ and such a hyperfunction $f(x)$ can be identified with the distribution $\wt f(x)$. When $n=1$, we only have to consider two cones $\Gamma_\pm=\pm(0,\infty)$ and we simply write \eqref{hyper-f} as $f(x)=F_+(x+i0)+F_-(x-i0)$. For example, the delta function and the Heaviside function
are given by $\delta(x)=(-2\pi i)^{-1}\big((x+i0)^{-1}-(x-i0)^{-1}\big)$ and 
$H(x)=(-2\pi i)^{-1}\big(\log(x+i0)-\log(x-i0)\big)$, where $\log z$ has slit along $(0,\infty)$.

We next define the singular part of hyperfunctions. We say that a hyperfunction $f(x)$
is {\em micro-analytic} at $(x_0;\im\xi_0)\in\im T^*\R^n\setminus\{0\}$ if $f(x)$ admits, near $x_0$, an expression  of the form \eqref{hyper-f} such that 
$\langle \xi_0, y\rangle<0$ for any $y\in\cup_j \Gamma_j$. The sheaf of {\em microfunctions} $\calC$ is defined as a sheaf on $\im T^*\R^n\setminus\{0\}$ with the
stalk at $(x_0;i\xi_0)$ given by the quotient space
$$
 \calC_{(x_0;i\xi_0)}=\calB_{x_0}/\{f\in\calB_{x_0}:
 \text{$f$ is micro-analytic at $(x_0;i\xi_0)$}\}.
$$
Since the definition of $\calC$ is given locally, we can also define the sheaf of microfunctions $\calC_M$ on $iT^*M\setminus\{0\}$ for a real analytic manifold $M$.

We now introduce a subclass of microfunctions that contains the Bergman and Szeg\"o kernels. Let $N\subset M$ be a real hypersurface with a real analytic defining function $\rho(x)$ and let  $Y$ be its complexification given by $\rho(z)=0$ in $X$.  Then, for each point $p\in N$, we consider a (multi-valued) holomorphic function of the form 
\begin{equation}\label{holo-micro}
 u(z)=\varphi(z)\rho(z)^{-m}+\psi(z)\log\rho(z),
\end{equation}
where $\varphi,\psi\in\calO_{X,p}$ and $m$ is a positive integer. A class modulo $\calO_{X,p}$ of $u(z)$ is called  a germ of {\em holomorphic microfunction} at $(p;i\xi)\in iT^*_N M\zs=\{(z;\lambda d\rho(z))\in T^*M:z\in N,\lambda\in\R \zs\}$, and we denote the sheaf of holomorphic microfunction on $iT^*_N M\zs$ by $\calC_{N|M}$. For a holomorphic microfunction $u$, we may assign a microfunction by taking the ``boundary values" from $\pm\Im \rho(z)>0$ with signature $\pm1$, respectively, as in the expression of $\delta(x)$ above,
which corresponds to $(-2\pi i\,z)^{-1}$. Thus we may regard $\calC_{N|M}$ as a subsheaf of $\calC_M$ supported on $i T^*_N M\zs$. With respect to local coordinates $(x',\rho)$ of $M$, each $u\in\calC_{N|M}$ admits a unique expansion
\begin{equation}\label{u-exp}
 u(x',\rho)=\sum_{j=k}^{-\infty} a_j(x')\Phi_j(\rho),
\end{equation}
where $a_j(x')$ are real analytic functions and
$$
\Phi_j(t)=
\begin{cases}j!\,t^{-j-1}\quad &\text{for }  j\ge0, 
\\ 
 \frac{(-1)^j}{(-j-1)!}\,t^{-j-1}\log t \quad &\text{for }  j<0.
\end{cases}
$$
If $u\ne0$ we may choose $k$ so that $a_k(x')\not\equiv0$
and call $k$ the {\em order of} $u$;
moreover, if $a_k(x')\ne0$ then we say that $u$ is 
{\em nondegenerate} at $(x',0)\in N$.

A differential operator
$P(x,D_x)=\sum a_\alpha(x)D^\alpha_x$, where $D^\alpha_x=(\pa/\pa x_1)^{\alpha_1}\cdots$ $(\pa/\pa x_n)^{\alpha_n}$, with real analytic coefficients acts on microfunctions; it is given by the application of the complexified operator $P(z,D_z)$ to each $F_j(z)$
in the expression \eqref{hyper-f}. 
Moreover, at $(p;i(1,0,\dots,0))\in \im T^*\R^n$, we can also define
the inverse operator $D_{x_1}^{-1}$ of $D_{x_1}$ by taking indefinite integrals of each $F_j$ in $z_1$. The microdifferential operators are defined as a ring generated by these operators. A germ of {\em microdifferential operator 
of order $m$} at $(x_0;i\xi_0)\in i T^*\R^n$ is a  series of
holomorphic functions 
 $\{P_j(z,\zeta)\}_{j=m}^{-\infty}$ defined on a conic neighborhood $U$ of
$(x_0;i\xi_0)$ in $T^*\C^n$ satisfying the following conditions:

(1) 
$P_j(z,\lambda\zeta)=\lambda^j P(z,\zeta)$ for
$\lambda\in\C\setminus\{0\}$;

(2)  For each compact set $K\subset U$, there
exists a constant $C_K>0$ such that
$\sup_K|P_{-j}(z,\zeta)|\le j!\,C_K^j$ for any $j\in\N=\{0,1,2,\dots\}$.

The series $\{P_j\}$ is denoted by $P(x,D_x)$,
and the formal series
$P(z,\zeta)=\sum P_j(z,\zeta)$ is called the {\em  total symbol},
while $\sigma_m(P)=P_m(z,\zeta)$ is called the {\em principal symbol}. 
The {\em product} and {\em adjoint} of microdifferential operators can be 
defined by the usual formula of symbol calculus:
$$
\begin{aligned}
 (PQ)(z,\zeta)&=\sum_{\alpha\in\N^n}\frac{1}{\alpha!}(D_\zeta^\alpha P(z,\zeta))D_z^\alpha Q(z,\zeta),\\
P^*(z,\zeta)&=\sum_{\alpha\in\N^n}
\frac{(-1)^{|\alpha|}}{\alpha!}
 D_z^\alpha D_\zeta^\alpha P(z,-\zeta).
\end{aligned}
$$
It is then shown that $P$ is invertible on a neighborhood of $(x_0;i\xi_0)$ if and only if $\sigma_m(P)(x_0,i\xi_0)\ne0$.

While these definitions based on the choice
of coordinates, we can introduce a transformation law of microdifferential operators
under coordinate changes and define the sheaf of the ring microdifferential
operators $\calE_M$ on $iT^*M$ for real analytic manifolds $M$. It  then turns out that
the adjoint depends only on the choice of volume element
$dx=dx_1\wedge\cdots\wedge dx_n$.

The action of differential operators on microfunctions can be extended
to the action of microdifferential operators so that
$\calC_M$ be a left $\calE_M$-module.
This is done by using the Laurent expansion of $P(z,\zeta)$ in $\zeta$
and then substituting $D_z$ and $D_{z_1}^{-1}$, or by introducing
a kernel function associated with the symbol (analogous to the
distribution kernel of a pseudodifferential operator).
Then $\calC_{N|M}$ becomes an $\calE_M$-submodule of $\calC_M$.
We can also define the right action of $\calE_M$ on $\calC_M\otimes\pi^{-1} v_M$, where $v_M$ is the sheaf of densities on $M$ and
$\pi\colon iT^*M\to M$ is the projection. It is given by
$(udx)P=(P^*u)dx$, where the adjoint is taken with respect to $dx$
(here $P^*$ depends on $dx$, but $(P^*u)dx$ is determined by $udx$).

We also consider microdifferential operators with a real analytic
parameter, that is, a $P=P(x,t,D_x,D_t)\in\calE_{M\times\R}$ 
that commutes with $t$. This is equivalent to saying that
the total symbol of $P$ is independent of the dual variable of
$t$; so we denote $P$ by  $P(x,t,D_x)$.
Note that $P(x,t,D_x)$, when $t$ is regarded as a parameter, acts on $\calC_M\otimes \pi^{-1}v_M$ from the right.

\subsection{Microfunctions associated with domains}
\label{microfunction-domain}
Now we go back to our original setting where $M$
is a complex manifold and $N=\partial\Omega$.
We have already seen that the Bergman kernel determines a
section of $\calC_{\pa\Omega|M}$, which we call the {\em local Bergman kernel} $B(x)$. Here $x$ indicate variables on
$M_\R$. Note that the local Bergman kernel
is defined for a germ of strictly pseudoconvex hypersurfaces.
Similarly, we can define the {\em local Szeg\"o kernel}\/: if we fix a real analytic surface element $d\sigma$ on $\pa\Omega$ 
and define the Szeg\"o kernel, then the
coefficients of the expansion \eqref{S-exp} are shown to be real analytic
and defines
a section $S(x)$ of $\calC_{\pa\Omega|M}$; see Remark \ref{Kas-rem} below.
We sometimes identify the surface element $d\sigma$ with the delta function
$\delta(\rho(x))$, or $\delta(\rho(x))dv$, normalized by $d\rho\wedge d\sigma=dv$. Note that the microfunction $\delta(\rho(x))$ corresponds to the holomorphic microfunction $(-2\pi\im\rho(z,\cw))^{-1}$
mod $\calO_X$, which we denote by $\delta[\rho]$.
Similarly, the Heaviside function $H(\rho(x))$ corresponds to 
a section $H[\rho]$ of $\calC_{\pa\Omega|M}$, 
which is represented by $(-2\pi\im)^{-1}\log\rho(z,\cw)$.

Our main object $\Vol(\Omega_\e)$ can be also seen as
a holomorphic microfunction.
In fact, since $u(\e)=\Vol(\Omega_\e)$ is a function
of the form $u(\e)=\varphi(\e)\e^{-n}+\psi(\e)\log\e$,
where $\varphi$ and $\psi$ are real analytic near $0$,
we may complexify $u(\e)$ and define a germ of holomorphic microfunction $u(\wt\e)\in\calC_{\{0\}|\R}$ at $(0;\im)$.
Note that $\Vol(\Omega_\e)\in \calC_{\{0\}|\R}$ is
expressed as an integral of the local Bergman kernel:
\begin{equation}\label{Vol-int} 
 \int B(x)H[\rho-\e](x)dv(x).
\end{equation}
Here $H[\rho-\e](x)$ is a section of $\calC_{\pa\wt\Omega|\wt M}$,
where
 $\wt\Omega=\{
(x,\e)\in\wt M= M\times\R:\rho(x)>\e\}$.
See Remark \ref{integral-rem} for the definition of this integral.

More generally, for a section $u(x,\e)$ of 
$\calC_{\pa\wt\Omega|\wt M}$ defined globally in $x$ for
small $\e$ and a global section $w(x)dx$ of
$\calC_{\pa\Omega|M}\otimes \pi^{-1}v_M$, we can define
 the integral of microfunction
$$
  \int u(x,\e)w(x)dx
$$
at $(0;\im )\in \im T^*\R$, which takes value in $\calC_{\{0\}|\R}$.
For such an integral, we have a formula of integration by part,
which is clear from the definition of the action of microdifferential operators in terms of kernel functions \cite{KKK}.

\begin{lem}\label{integral-lemma}
If $P(x,\e,D_x)$ is a microdifferential operator
defined on a neighborhood of the support of $u(x,\e)$,
then
\begin{equation}\label{integral-by-part}
 \int\big(P u\big)wdx= \int u\big(wdx\,P\big).
\end{equation}
\end{lem}

\begin{rem}\rm\label{integral-rem}
We here recall the definition of the integral \eqref{Vol-int} and
show that it agrees with $\Vol(\Omega_\e)$.
For a general definition of the integral of microfunctions, we refer to \cite{KKK}. 
 Write $dv=\lambda  d\rho\wedge d\sigma$
 and complexify $\lambda(x',\rho)$ to $\lambda(x',\wt\rho)$ for $\wt\rho\in\C$
 near $0$.
Then, define a holomorphic function $f(\wt\e)$ on $\Im\wt\e>0$, 
$|\wt\e|\ll1$,
by the path integral 
\begin{equation}\label{path-int}
 f(\wt\e)=\int_{\pa\Omega}\int_{\gamma_1} 
 B(x',\wt\rho)\frac{1}{2\pi i}\log(\wt\rho-\wt\e)
 \lambda(x',\wt\rho)\,d\wt\rho d\sigma(x'),
\end{equation} 
where $\gamma_1$ is a path connecting $a$ and $b$, with $a<0<b$, 
such that the image is contained in $0<\Im\wt\rho<\Im\wt\e$
except for the both ends.
Then \eqref{Vol-int} is given by
 $f(\e+\im0)\in\calC_{\R},{}_{(0;i)}$, which
 is independent of the choice of $a,b$ and $\gamma$.
We now show $f(\e+i0)=\Vol(\Omega_\e)$ as a microfunction.
For each $\wt\e$ with $\Im\wt\e>0$,
choose another path
connecting $b$ and $a$ so that $\gamma_2\gamma_1$ be
a closed path surrounding $\wt\e$ in the positive direction.
Since the integral along $\gamma_2$ gives a function
that can be analytically continued to $0$, we may replace
$\gamma_1$ in \eqref{path-int} by $\gamma_2\gamma_1$
without changing its $\calO_{\C,0}$ modulo class.
Now restricting $\wt\e$ to positive real axis,
and letting the path $\gamma_2\gamma_1$ shrink to the line
segment $[\e,b]$, we see that
$f(\e)$ agrees with $\Vol(\Omega_\e)$ modulo analytic functions
at $0$.
\end{rem}

\subsection{Quantized contact transformations}
\label{quantized-transf}
We finally recall a property of holomorphic microfunctions that follows from the strictly pseudoconvexity of $\pa\Omega$. Let $z$ be local holomorphic coordinates of $M$. Then we write $P(x,D_x)=P(z,D_z)$ (resp.\ $P(\cz,D_\cz)$) if $P$ commutes with $\cz_j$ and $D_{\cz_j}$ (resp.\ $z_j$ and $D_{z_j}$). Similarly for $P(x,t,D_x,D_t)\in\calE_{M\times\R}$ we write, e.g.\ $P(x,t,D_x,D_t)=P(z,t,D_z)$ if $P$ commutes with $\cz_j$, $D_{\cz_j}$ and $t$. Clearly, the class of operators
$P(z,D_z)$ and $P(\cz,D_\cz)$ are determined by the complex structure of $M$.

\begin{lem} 
Let $N$ be a strictly pseudoconvex hypersurface in $M$ with a defining function $\rho$.  Then for each section $u$ of\/ $\calC_{N|M}$, there exists a unique microdifferential operator $R(z,D_z)$ such that $u=R(z,D_z)\delta[\rho]$. Moreover,  $u$ and $R$ have the same order, and $u$ is nondegenerate if and only if $R$ is invertible.
\end{lem} 

Note that the same lemma holds when $\delta[\rho]$ is replaced by $H[\rho]$, or more generally, by a nondegenerate section $u$ of $\calC_{N|M}$, except for the statement about the order.

The strictly pseudoconvexity of $N$ implies that the projection $p_1:T_Y^*X\subset T^*(M\times\conj M)\to T^*M$ is a local biholomorphic map, where $Y$ is the complexification of $N$ in $X$.  The surjectivity and the injectivity of $p_1$ imply the existence and uniqueness of $R(z,D_z)$, respectively. If we apply the same argument for $\cz$, we obtain a local biholomorphic map $p_2:T_Y^*X\to T^*\conj M$ and a contact (or homogeneous symplectic) transformation
$\phi(z,\zeta)=p_1\circ p_2^{-1}(z,-\zeta)$.
Then, for each nondegenerate $u$, the lemma above gives a map $\Phi:P(z,D_z)\mapsto Q(\cz,D_\cz)$
such that $(P-Q^*)u=0$, where the adjoint is taken with respect to $|dz|^2$.  This is an isomorphism of the rings
and satisfies $\sigma_m(\Phi(P))\circ\phi=\sigma_m(P)$ if $P$ has order $m$; hence $\Phi$ is called a {\em quantized contact transformation}  with a {\em generating function} $u$.
It is shown that a quantizations of $\phi$ determines a generating function uniquely up to a constant multiple.
Chapter 1 of \cite{Schapira} is a good reference for this subject.

\section{Kashiwara's analysis of the kernel functions}  
\label{sec-Kashiwara}
In this section we recall Kashiwara's analysis of the Bergman kernel and its analogy for the Szeg\"o kernel. Then we derive some microdifferential equations satisfied by these kernels.

\subsection{A relation between the local Bergman and Szeg\"o kernels}
Under the formulation of the previous section, 
Kashiwara's theorem \cite{Kas} for the Bergman kernel and its analogy for the Szeg\"o kernel can be stated as follows:

\begin{thm}\label{kashiwara-thm}
{\rm (i)} 
The local Bergman kernel satisfies
$$
\big(P(z,D_z)-Q(\cz,D_\cz)\big)B=0
$$
for any pair of microdifferential operators $P(z,D_z)$ and $Q(\cz,D_\cz)$ such that
\begin{equation}\label{q-trB}
 (H[\rho]dv)(P(z,D_z)-Q(\cz,D_\cz))=0.
\end{equation}
Moreover, the local Bergman kernel is uniquely determined by this property up to a constant multiple.

{\rm (ii)} 
The local Szeg\"o kernel satisfies
$$
\big(P(z,D_z)-Q(\cz,D_\cz)\big)S=0
$$
for any pair of microdifferential operators $P(z,D_z)$ and $Q(\cz,D_\cz)$ such that

\begin{equation}\label{q-trS}
(\delta[\rho]dv)(P(z,D_z)-Q(\cz,D_\cz))=0.
\end{equation}
Moreover, the local Szeg\"o kernel is uniquely determined by this property up to a constant multiple.
\end{thm}

\begin{rem}\rm\label{Kas-rem}
In \cite{Kas},  Kashiwara stated (i) and gave its heuristic proof, which can be equally applied to (ii). Also, as a premise for this theorem, he stated the real analyticity of the coefficients of the asymptotic expansion of the Bergman kernel, while the proof has not been published. Now a proof of these theorem and claim, based on Kashiwara's lectures, is available in Kaneko's lecture notes \cite{Kan}; the arguments there can be also applied to the case of the Szeg\"o kernel.
\end{rem}

Take holomorphic coordinates $z$ and
write $dv=\varphi |dz|^2$. Then \eqref{q-trB} can be rewritten as
$$
 (P^*-Q^*)\varphi H[\rho]=0,
$$
where the adjoint is taken with respect to $|dz|^2$. It follows that the maps $P^*(z,D_z)\mapsto Q(\cz,D_\cz)$ and $Q(\cz,D_\cz)\mapsto P^*(z,D_z)$ are the quantized contact transformations generated by $\varphi H[\rho]$ and $B$,
respectively, and are the inverse of each other. Thus we can say that the theorem states the reproducing property of the kernel on microdifferential operators. In particular, we see that the uniqueness statement of the theorem follow from that of the generating function.

From Theorem \ref{kashiwara-thm}, we can easily derive a microdifferential relation between the local Bergman and
Szeg\"o kernels.

\begin{prop}\label{B-S-relation}
 Let $R(z,D_z)$ be a microdifferential operator such 
that $(H[\rho]dv) R=\delta[\rho]dv$.
Then $RS=B$.
\end{prop}

{\em Proof.}
We first show that $(P(z,D_z)-Q(\cz,D_\cz))RS=0$
for any pair $P$ and $Q$ satisfying $(H[\rho]dv)(P-Q)=0$.
Noting that $Q(\cz,D_\cz)$ commutes with $R(z,D_z)$, we see
from
$H[\rho]dv=(\delta[\rho]dv) R^{-1}$ that
$
(\delta[\rho]dv)(R^{-1}PR-Q)=0$.
Since $R^{-1}PR$ is an operator of $z$-variable,
Theorem \ref{kashiwara-thm} implies 
$(R^{-1}PR-Q)S=0$ and thus $(P-Q)RS=0$. 
Now by the uniqueness statement
of Theorem \ref{kashiwara-thm}, we have $B=c\,RS$ for 
a constant $c$.

It remains to show that $c=1$.
This can be done by computing explicitly the leading term of these
kernels. Take local coordinates $z=(z',z_n)$ such that
the boundary $\pa\Omega$ is locally given
by the defining function
$$
 \rho_0(z,\cz)=z_n+\cz_n-z'\cdot\cz'+F(z,\cz), 
\quad
F=O(|z|^3).
$$
Then write $\rho=e^{-f(z,\cz)}\rho_0$ and $dv=e^{g(z,\cz)}dV$, where
$dV$ is the standard measure on $\C^n$.
Since $\Omega$ is osculated at $0$ to the third order by the 
Siegel domain, we see that
$$
B=\frac{n!}{\pi^n}e^{-g}\rho_0^{-n-1}\big(1+O(1)\big)
\qtext{and}
S=\frac{(n-1)!}{\pi^n}e^{-f-g}\rho_0^{-n}\big(1+O(1)\big),
$$
where $O(1)$ denotes a term that vanishes at $z=0$;
see \cite{BS}.
On the other hand, setting $f_0=f(0,0)$ and $g_0=g(0,0)$, we
have $e^{f_0}D_{z_n}e^gH[\rho]=e^g\delta[\rho]+u$ for a
degenerate germ $u\in\calC_{\pa\Omega|M}$ at $(0;id\rho)$.
Thus $R(z,D_z)=-e^{f_0}D_{z_n}
+P(z,D_z)$, where $P$ has order at most $1$ and $\sigma_1(P)(z,\zeta)$
vanishes at $(0;id\rho)$.
Using the expression of $S$ above, we have
$$
RS=-e^{f_0}D_{z_n}S+PS=\frac{n!}{\pi^n}e^{-g}\rho_0^{-n-1}\big(1+O(1)\big),
$$
which implies $c=1$.
\qed

\subsection{Variational formula of the local Szeg\"o kernel}
 Let $\{\Omega_t\}_{t\in I}$ be a $C^\omega$ family of strictly pseudoconvex domains with $C^\omega$ boundaries,
where $I\subset\R$ is an open interval. 
Here a $C^\omega$ family means that 
$\wt\Omega=\{(x,t)\in \wt M=M\times I:x\in\Omega_t\}$
admits a $C^\omega$ defining function $\rho_t(x)$ such that $d_x\rho_t(x)\ne0$ on $\pa\wt\Omega$. If we fix $\rho_t$,
we can assign for each $\pa\Omega_t$ a surface element
$d\sigma_t$ by $\delta[\rho_t]dv$.
We here consider the microdifferential equations for the family of the local Szeg\"o kernels of $(\pa\Omega_t,d\sigma_t)$.

\begin{prop}\label{Dt-prop} 
There exists a section $S_t(x)$ of $\calC_{\pa\wt\Omega|\wt M}$
such that, for each $t$, $S_t(x)$  gives the local Szeg\"o kernel of $(\pa\Omega_t,d\sigma_t)$.
Moreover, $S_t(x)$ satisfies
\begin{equation}\label{St-system}
\big(P(z,t,D_z)-Q(\cz,D_\cz,D_t)\big)S_t(x)=0
\end{equation}
for any pair of microdifferential operators $P(z,t,D_z)$ and $Q(\cz,D_\cz,D_t)$ such that
$$
(\delta[\rho_t]d\wt v)\big(P(z,t,D_z)-Q(\cz,D_\cz,D_t)\big)=0,
$$
where $d\wt v=dv\wedge dt$ is a volume element on $\wt M$.
In particular, if $R(z,t,D_z)$ satisfies $D_t  \delta[\rho_t]dv=(\delta[\rho_t]dv) R(z,t,D_z)$, then
\begin{equation}\label{DtS}
-D_tS_t= R(z,t,D_z)S_t.
\end{equation}
\end{prop} 

Analogous proposition for the local Bergman kernel has been given in
\cite{Hi1}, where we considered the family of local Bergman kernels $B_t(x)$ of
$(\Omega_t,|dz|^2)$ for domains in $\C^n$ and obtained the exactly same statement for $B_t$ with $H[\rho_t]|dz|^2$ in place of $\delta[\rho_t]dv$. In particular, we have
\begin{equation}\label{DtB}
-D_tB_t=\wt R(z,t,D_z)B_t
\end{equation}
for $\wt R$ satisfying $D_t H[\rho_t]|dz|^2=(H[\rho_t]|dz|^2) \wt R(z,t,D_z)$.
We here use Proposition \ref{B-S-relation}
to translate this formula into the one for $S_t$.

\medskip

{\em Proof.} 
First note that the ring of operators of the form
$Q(\cz,D_\cz,D_t)$ is the ring generated by $\cz_1,\dots,\cz_n,D_{\cz_1},\dots,D_{\cz_n},D_t$, and hence                                                                                                                        it suffices to prove \eqref{St-system}
when $Q$ is one of these generators. For $\cz_j$ and $D_{\cz_j}$,
it is clear from Theorem \ref{kashiwara-thm}.
To prove the case $Q=D_t$,
take $A(z,t,D_z)$ such that $\delta[\rho_t]dv=(H[\rho_t]dv)A$
and compute
$$
\begin{aligned}
 (\delta[\rho_t]d\wt v)D_t
  &=(H[\rho_t]d\wt v)AD_t
\\
&= (H[\rho_t]d\wt v)[A,D_t]+(H[\rho_t]d\wt v)D_tA\\
  &=(H[\rho_t]d\wt v)([A,D_t]-\wt RA)\\
  &=(\delta[\rho_t]d\wt v)A^{-1}([A,D_t]-\wt RA).
\end{aligned}
$$
Since $[t,[D_t,A]]=0$, 
we have 
$R(z,t,D_z)=A^{-1}(\wt RA-[A,D_t])$.
On the other hand, Proposition \ref{B-S-relation} implies
$B_t=AS_t$ and thus
$$
\begin{aligned}
 R S_t&=A^{-1}(\wt RA-[A,D_t])S_t\\
  &=(A^{-1}\wt RA+A^{-1}D_tA-D_t)S_t\\
  &=A^{-1}(\wt R+D_t)B_t-D_tS_t.
\end{aligned}
$$
Therefore, using \eqref{DtB}, we get $RS_t=-D_tS_t$.
\qed

\section{Proofs of the main theorems}   
\label{sec-proofs}
Now we are ready to prove the main theorems.
We first note that the theorems can be reduced to the ones in the real analytic category by approximations.  The key fact is that
the asymptotic expansion up to each fixed order of the Bergman and Szeg\"o kernels are determined by the finite jets of 
$\rho$, $d\sigma$ and $dv$ at each boundary point.
Thus, for a domain $\Omega$ with $C^\infty$ defining function $\rho$
and the contact form $\th=\im(\pa\rho-\conj\pa\rho)$ on $\pa\Omega$,
by taking a series of $C^\omega$ functions $\rho_{j}$ that converges to $\rho$ in $C^k$-norm for any $k$, we may express $L(\pa\Omega,\th)$ as the limit of $L(\pa\Omega_j,\th_j)$, where $\Omega_j=\{\rho_j>0\}$. 
To reduce Theorem \ref{thmSzego} (i)
to the real analytic case, we only have to take another sequence
of $C^\omega$ contact forms $\{e^{f_j}\th_j\}$ approximating a
given contact form $e^f\th$ so that $L(\pa\Omega_j,e^{f_j}\th_j)=L(\pa\Omega_j,\th_j)$ implies $L(\pa\Omega,e^{f}\th)=L(\pa\Omega,\th)$.
Similar arguments of approximation can be applied to the other cases. 

In the following we prove the theorems in the real analytic category.

\subsection{Proof of Theorem \ref{thmSzego}}
Taking a $C^\omega$ family of defining functions $\rho_t(x)$,
we define 
$S_t(x)$ to be the local Szeg\"o kernel for the surface element given by  $\delta[\rho_t]dv$.
Let $\wt\Omega=\{(x,\e,t)\in \wt M=M\times\R^2:\rho_t(x)>\e\}$ so that
$\delta[\rho_t-\e]$ defines a section of $\calC_{\pa\wt\Omega|\wt M}$ and
consider the integral 
$$
 A(\e,t)=\int S_t(x)\delta[\rho_t-\e](x)dv(x),
$$
which is well-defined as a germ of $\calC_{\{0\}\times\R|\R^2}$ at $(0,0;\im(1,0))$. Write
$$
A(\e,t)=\varphi(\e,t)\e^{-n}+\psi(\e,t)\log\e
$$
and set $L_t=\psi(0,t)$, which we call the coefficient of
$\e^0\log\e$. Then
our goal is to 
prove the independence of $L_t$ on $t$, because it contains the theorem:
For the statement (i), we take $\rho_t=e^{tf}\rho$ so that $\delta[\rho_0]dv$ and
$\delta[\rho_1]dv$ correspond to $\th\wedge (d\th)^{n-1}$ and $\wt\th\wedge (d\wt\th)^{n-1}$ respectively;
then $L_0=L(\pa\Omega,\th)$ and $L_1=L(\pa\Omega,\wt\th)$ agree.
For the statement (ii), we have $L_t=L(\pa\Omega_t)$, which is independent of $t$.

Take a microdifferential operator $R(z,t,D_z)$ such that
$$
 D_t\delta[\rho_t]dv=(\delta[\rho_t]dv)R.
$$
Then Proposition \ref{Dt-prop} implies $-D_tS_t=RS_t$.
Using this and \eqref{integral-by-part}, we have
$$
\begin{aligned}
 D_t A(t,\e)&=D_t\int S_t\,\delta[\rho_t-\e]dv\\
 &=\int (D_tS_t)\delta[\rho_t-\e]+S_tD_t\delta[\rho_t-\e]dv\\
&=\int (-RS_t)\delta[\rho_t-\e]+S_tD_t\delta[\rho_t-\e]dv\\
&=\int S_t\big(-(\delta[\rho_t-\e]dv) R+D_t\delta[\rho_t-\e]dv\big).
\end{aligned}
$$ 
Since $D_t\delta[\rho_t-\e]dv-(\delta[\rho_t-\e]dv) R$ vanishes
at $\e=0$, we may take, by using Lemma \ref{Division-lemma} below, a section $u$ of $\calC_{\pa\wt\Omega|\wt M}$ such that 
$$
 D_t\delta[\rho_t-\e]dv-(\delta[\rho_t-\e]dv) R=\e\,u dv.
$$
Hence we have
$$
 D_t A(t,\e)=\e\int\! S_t(x)\, u(x,t,\e)dv(x).
$$
The integral in the right-hand side takes value in
$\calC_{\{0\}\times\R|\R^2}$;
thus the right-hand side does not contain $\e^0\log\e$
term. This
implies $D_t L_t=0$.
\qed

\begin{lem}\label{Division-lemma}
Let $\{\Omega_\e\}$ be a real analytic family of domains
in $M$ and $\rho(x,\e)$ be the defining function of
$\wt\Omega=\{(x,\e)\in M\times\R:x\in\Omega_\e\}$ such that
$d_x\rho(x,\e)\ne0$ on the boundary.
If $u(x,\e)\in\calC_{\pa\wt\Omega|\wt M}$ satisfies $u(x,0)=0$ in $\calC_{\pa\Omega_0|M}$, 
then there exists a germ $v\in\calC_{\pa\wt\Omega|\wt M}$ such that $u=\e\,v$. 
\end{lem}

{\em Proof.} Take coordinates $(x',\rho,\e)$ for $\wt M$
and expand $u$ as in \eqref{u-exp} with coefficients $a_j(x',\e)$.
Then $u(x',\rho,0)=0$ implies $a_j(x',0)=0$ so that $a_j=\e a_j'$ for real analytic functions $a_j'(x,\e)$.
Thus we may set $v=\sum a_j'(x',\e)\Phi_{j}(\rho)$. 
\qed

\subsection{Proof of Theorem \ref{thmBergman}}
Take a microdifferential operator $R(z,\e,D_z)$ such that
\begin{equation}\label{HR-exp}
(H[\rho-\e]dv) R=\delta[\rho-\e]dv.
\end{equation}
Then $(H[\rho]dv) R(z,0,D_z)=\delta[\rho]dv$
and hence $R(z,0,D_z)S=B$ by Proposition \ref{B-S-relation}.
So, applying Lemma \ref{Division-lemma} for $\Omega_\e=\Omega$, we have
\begin{equation}\label{RS-exp}
 R(z,\e,D_z) S(x)=B(x)+\e\, B'(x,\e),
\end{equation}
where $B'$ for a section of $\calC_{\pa\Omega\times\R|M\times \R}$.
Using \eqref{HR-exp} and \eqref{RS-exp}, we compute
$$
\begin{aligned}
\int S\delta[\rho-\e]dv
&=\int S\big(( H[\rho-\e]dv) R\big)\\
&=\int (R S) H[\rho-\e]dv\\
&=\Vol(\Omega_\e)+\e \int B'H[\rho-\e]dv.
\end{aligned}
$$
Since the integral in the right-hand side takes value in
$\calC_{\{0\}|\R}$, its $\e$ multiple cannot contain 
$\e^0\log\e$ term.
Therefore the coefficients of $\e^0\log \e$ of $\int S\,\delta[\rho-\e]dv$ and $\Vol(\Omega_\e)$ agree; the former gives $L(\pa\Omega)$
and the theorem follows.
\qed

\medskip

\section{Proof of Proposition \ref{prop-tube-volume}} 
\label{sec-tube} 

Let $\omega$ be the curvature of $E$, which is assumed to give a K\"ahler form, and set $dv_X=\omega^{n-1}/(n-1)!$.
Then $\pi^*\omega=-i\pa\conj\pa\rho$ and  thus $dv=d\rho\wedge d\phi\wedge  \pi^* dv_X$,
where $\pi\colon E^*\to X$ is the projection and $\phi=\arg z_0$ for a fiber coordinate $z_0$ of $E^*$.
From this formula, we see that the surface element $d\sigma=d\phi\wedge\pi^* dv_X$ on $\pa\Omega$ corresponds to 
$\delta[\rho]dv$.

 Let $\calA_m(\pa\Omega)$ be the subspace of
  $\calA(\pa\Omega)=\calA(\pa\Omega,d\sigma)$ consisting of 
  functions that are homogeneous of degree $m$,
 that is,  $\varphi(\lambda v)=\lambda^m\varphi(v)$ for any $\lambda\in\C$ with $|\lambda|<1$. 
Then we have an orthogonal decomposition
$
\calA(\pa\Omega)=\bigoplus_{m=0}^\infty\calA_m(\pa\Omega)
$.
Here each $\calA_m(\pa\Omega)$ can be canonically identified with
$H^0(X,E^{\otimes m})$ and hence has finite dimension $d_m$.
Thus taking, for each $m$,  an orthonormal
basis $\varphi_{1,m},\dots,\varphi_{d_m,m}$ of $\calA_m(\pa\Omega)$,
we may form
 a complete orthonormal system $\{\varphi_{j,m}\}_{j,m}$ of 
$\calA(\pa\Omega)$. 
Since $dv$ has singularity along the zero section, the
constant function $\varphi_{1,0}$ is not contained
in $L^2(\Omega,dv)$; but except for that,
all $\varphi_{j,m}$, $m>0$, are also contained in $\calA(\Omega)=\calA(\Omega,dv)$
and give a complete orthogonal basis of $\calA(\Omega)$. Therefore
 $\|\varphi_{j,m}\|^2_{\calA(\Omega)}=1/m$ implies
\begin{equation}\label{SandB}
 S(v,\conj v)=\sum_{j,m}|\varphi_{j,m}(v)|^2,
 \quad
 B(v,\conj v)=\sum_{j,m} m|\varphi_{j,m}(v)|^2.
\end{equation}
Since $B(v,\conj v)$ and $S(v,\conj v)$ are invariant under the $S^1$-action, we may write 
 $B(v,\cv)=B_x(\rho)$ and $S(v,\cv)=S_x(\rho)$,
 where $\pi(v)=x$ and $\rho=\rho(v)$.
In the coordinates $(x,\rho,\phi)$, we have $D_\rho |\varphi_{j,m}|^2=-m |\varphi_{j,m}|^2$ so that
\begin{equation}\label{DS=B}
-D_\rho S_x(\rho)=B_x(\rho).
\end{equation}
Thus integrating $B_x(\rho)dv$ on each fiber of $\Omega_\e\to X$, we have
\begin{equation}\label{B2S}
\Vol(\Omega_\e)=2\pi \int_X  
\Big(S_x(\e)-S_x(\infty)\Big) dv_X. 
\end{equation}
Recalling \eqref{u-exp}, we expand $S_x(\rho)$ as
$$
 S_x(\rho)\sim
 \sum_{j=n-1}^{-\infty} a_j(x)\Phi_{j}(\rho),
 \quad\text{where }a_j\in C^\infty(X),
$$
or equivalently, we write $S_\rho(x)$ as the Laplace transform
\begin{equation}\label{S-laplace}
 S_x(\rho)=\int_0^\infty e^{-t\rho}a(x,t)dt
\end{equation}
of a classical symbol $a(x,t)\in S^n(X\times\R_+)$ 
with asymptotic expansion $a(x,t)\sim\sum_{j=n-1}^{-\infty} a_j(x)t^{j}$
 at $t=\infty$.
Substituting this into \eqref{B2S} gives
$$
 \Vol(\Omega_\e)\sim2\pi\int_0^\infty \int_X e^{-t\e}a(x,t)dv_X dt
 \quad \text{as }\e\to 0.
$$

It only remains to prove 
\begin{equation}\label{a2P}
  \int_X a(x,m)dv_X\sim 
  d_m
  \quad\text{as }
 m\to\infty.
\end{equation}
Note that the orthogonal projection $\calA(\pa\Omega)\to\calA_m(\pa\Omega)$ is given by $f(v)\mapsto(2\pi)^{-1} \int_{S^1}e^{-im\phi}f(e^{i\phi}v)d\phi$. Thus we have
$$
  \sum_{j=1}^{d_m}|\varphi_{j,m}(v)|
  =\frac{1}{2\pi}\int_0^{2\pi} e^{-\im m\phi}S(e^{i\phi}v,\cv)d\phi.
 $$
We restrict this formula to $\pa\Omega$.
Then the right-hand side depends only on $x$ and define
a function $B_m(x)$, the Bergman kernel of $H^0(X,E^{\otimes m})$.
On the other hand, we can compute
the singularity of $S(e^{i\phi}v,\cv)$ at $\phi=0$ by complexifying the expression \eqref{S-laplace}: for small $\phi$, 
we have $\rho(e^{i\phi}v,\cv)=-i\phi$ so that
$$
 S(e^{i\phi}v,\cv)=\int_0^\infty e^{it\phi}a(x,t)dt.
$$
It is now clear from Fourier's inversion formula that 
\begin{equation}\label{BsimA}
B_m(x)\sim a(x,m)\quad\text{ as }
 m\to\infty.
\end{equation}
Since $\int_X B_m(x)dv_X=d_m$, we get \eqref{a2P}

\begin{rem} \rm \label{catlin-volume}
Our choice of volume element $dv$ is different form that used
in  Catlin \cite{Cat}, where he employed a volume element 
$e^{-\rho}dv$ which is also smooth along the zero section.
By following the argument deriving \eqref{SandB}, we can easily show that Catlin's Bergman kernel is given by the sum 
$B+S$ of our Bergman and Szeg\"o kernels.
 
\end{rem}

\begin{rem} \label{rem-volume}\rm
The relation \eqref{DS=B}, as microfunctions, can be also derived 
from Proposition \ref{kashiwara-thm}. In fact,
noting the fact the $\rho-2i\phi$ is holomorphic, we take
$P=-D_\rho-(i/2) D_\phi$; then $\delta[\rho]dv=H[\rho]dvP$, which implies \eqref{DS=B} since $D_\phi S(v,\cv)=0$.
Note also that the boundary singularity of $B(v,\cv)$
can be localized.  Thus we can modify $dv$ near the zero section
and make it smooth without changing \eqref{DS=B} modulo
smooth error.

\end{rem}

\end{document}